\documentclass[12 pt,a4paper]{amsart} 
\usepackage{amsfonts,amssymb,amscd,amsmath,enumerate,verbatim,calc} 
\usepackage{float}
\usepackage{amsthm}
\newtheorem{theorem}{Theorem}[section]
\newtheorem{definition}{Definition}[section]

\newtheorem{lemma}[theorem]{Lemma}
\newtheorem{corollary}[theorem]{Corollary}
\newtheorem{proposition}[theorem]{Proposition}
\newtheorem{remark}[theorem]{Remark}
\usepackage{tikz-cd}
\usepackage{mathtools}
\usepackage{xfrac}
\usepackage{csquotes}
\usepackage[all,cmtip]{xy}
\newtheorem{examples}[theorem]{Examples}

\numberwithin{equation}{section}

\usepackage{amsfonts}

\newcommand{\im}{\rm Im}

\renewcommand{\ker}{\textnormal{ker}}
\renewcommand{\im}{\textnormal{Im}}

\begin{document}
	\title{Relative Lie central extension and Schur multiplier of pairs of multiplicative Lie algebras}
	\author{Dev Karan Singh$^{1}$, Shiv Datt Kumar$^{2}$\vspace{.4cm}\\
		{$^{1}$Department of Mathematics, \\ Institute of Integrated Learning in Management University,\\ Greater Noida (UP) - 201306, India}\\
		{$^{2}$Department of Mathematics, \\ Motilal Nehru National Institute of Technology Allahabad \\ Prayagraj (UP) - 211004, India}}
	
	\thanks{$^1$devkaransingh1811@gmail.com, $^2$sdt@mnnit.ac.in.  }
	
	\begin{abstract}
		In this paper, we introduce the concept of relative Lie central extension for pair of multiplicative Lie algebras. Then, we discuss the concept of isoclinism for relative Lie central extensions and prove some related results. We also define the Frattini subalgebra for multiplicative Lie algebras and discuss its properties, finally the Schur multiplier for pair of multiplicative Lie algebras is introduced and under certain conditions prove the existence of multiplicative covering pair.
	\end{abstract}
	\maketitle
	\section{Introduction}
	
	 In 1991, G. J. Ellis, introduced the notion of multiplicative Lie algebra which is a generalization of groups and Lie algebras or Lie rings, to prove that five universal commutator identities generate all universal identities between commutators of weight $n$. 
	In \cite{FPW}, Point and Wantiez discussed the nilpotency and solvability of multiplicative Lie algebras. In \cite{AGNM}, Bak et al. developed two homology theories for multiplicative Lie algebras and explored their application to K-theory. In \cite{RLS}, Lal and Upadhyay established Schur Hopf's formula and defined the Lie exterior square of multiplicative Lie algebras. In \cite{GNM}, Donadze et al. introduced the non-abelian tensor product of these algebras and proved an analogue of Miller's theorem for groups.
	In \cite{AMS}, Kumar et al. discussed the capability and covers of multiplicative Lie algebras. In \cite{DMS}, Singh et al. covered the central extension theory of perfect multiplicative Lie algebras and local nilpotency, while in \cite{MS}, Pandey and Upadhyay defined the concept of isoclinism in multiplicative Lie algebras. In \cite{DASKumar} Singh et al. extend the concept of isoclinism for central extensions of multiplicative Lie algebras. In \cite{DKSD}, Singh and Kumar discuss the inducibility problem for automorphisms of multiplicative Lie algebra extensions and establish the Wells type exact sequences.
	
	
	In \cite{Loday}, Loday, introduced the concept of relative central extension of pairs of groups, which generalizes the classical notion of central extensions. In \cite{Ellis_Pair}, Ellis defined the Schur multiplier of a pair of groups, which is helpful for studying pairs of groups and third integral homology. In \cite{A. Pour}, \cite{chitti}, authors, discussed the properties of Schur multipliers and the covers of pairs of groups. In \cite{Arabyani}, Arabyani and Mohseni, introduced the concept of isoclinism for relative central extension of pairs of Lie algebras. In \cite{Safa}, Safa and Arabyani, discussed the c-covering pairs and some results on the existence of c-covers of a pair of Lie algebras.
	
	In this paper, we introduced the concept of relative Lie central extensions for pairs of multiplicative Lie algebras. This concept is the generalization of relative central extension of groups and Lie algebras, since action of multiplicative Lie algebras involves the notion of actions of both. We also discussed some properties of isoclinism for relative Lie central extension and Schur multiplier of a pair of multiplicative Lie algebras.
	
	This paper is organised as follows: In section 2, we recall some basic definitions of multiplicative Lie algebras and properties of action of multiplicative Lie algebras. In Section 3, we introduce the concept of the relative Lie central extension for a pair of multiplicative Lie algebras and discuss the properties of the Frattini subalgebra in the context of multiplicative Lie algebras. In Section 4, we define the concept of isoclinism for relative Lie central extensions of multiplicative Lie algebras, construct a relative multiplicative central extension from two given relative Lie central extensions, and obtain different types of isoclinic morphisms. In Section 5, we first define the Schur multiplier of pairs of multiplicative Lie algebras and then prove a four-term exact sequence relating the Schur multipliers of such pairs. We also define the covering pair of multiplicative Lie algebras and show that it is the homomorphic image of a quotient of a free multiplicative Lie algebra. Furthermore, we establish the existence of a covering pair under certain conditions.

	\section{Preliminaries}
	In this section, we recall definitions and results related to multiplicative Lie algebras.
	\begin{definition}\cite{GJ}
		A multiplicative Lie algebra is a triple $(G, \cdot, \star )$, where $(G,\cdot)$ is a group and binary operation $\star$ satisfy the following identities
		\begin{enumerate}
			\item $ g\star g=1 $  
			\item $ g\star(hk) = (g\star h){^h(g\star k)}$ 
			\item $ (gh)\star k = {^g(h\star k)} (g\star k)$ 
			\item $ ((g\star h)\star {^hk})((h\star k)\star{^kg})((k\star g)\star{^gh})=1$ 
			\item $ ^k(g\star h)=(^kg\star {^kh})$ 
		\end{enumerate}
		for all $g,h,k\in G$, where $^gh$ denotes $ghg^{-1}$ and $\star$ is called Lie product. 
	\end{definition}
	\noindent \textbf{Note:} If $(G,\cdot)$ is an abelian group, then $(G,\cdot,\star)$ is called a Lie ring.

	Let $(G, \cdot, \star)$ be a multiplicative Lie algebra. A subgroup $H$ of $G$ is called a \emph{subalgebra} if it is closed under the Lie product $\star$. If, in addition, $H$ is normal in $G$ and satisfies $g \star h \in H$ for all $g \in G$ and $h \in H$, then $H$ is called an \emph{ideal}. The ideal generated by all Lie products $g \star g'$ for $g, g' \in G$ is denoted by $G \star G$. A group homomorphism $\psi : G \to G'$ between multiplicative Lie algebras is said to be a \emph{multiplicative Lie algebra homomorphism} if it preserves the Lie product.
	The group center $Z(G) = \{ g \in G \mid [g, g'] = 1 \text{ for all } g' \in G \}$ and the Lie center $LZ(G) = { g \in G \mid g \star g' = 1 \text{ for all } g' \in G }$ are both ideals of $G$, and their intersection $\mathcal{Z}(G) = Z(G) \cap LZ(G)$ is also an ideal and called \emph{multiplicative center} of $G$. The ideal $(G\star G)[G,G]$ generated by all elements of the form $[g_1, g_2](g_3 \star g_4)$, where $g_i \in G$, is denoted by ${}^M[G, G]$ and called \emph{multiplicative commutator}. The \emph{lower central series} $\{M_n(G)\}$ is defined inductively by $M_0(G) = G$, $M_1(G) = {^M[G, M_0(G)]}$ and $M_{n+1}(G) = [G, M_n(G)]$. The multiplicative Lie algebra $G$ is said to be \emph{nilpotent} if there exists some $n \in \mathbb{N}$ such that $M_n(G) = {1}$.

	\begin{definition} \cite{GNM}
		Let $G$ and $H$ be two multiplicative Lie algebras. By an action of $G$ on $H$ we mean an underlying group action of  $G$ on  $H$ and $H$ on  $G$, together with a map $G \times H \to H, (g, h) \mapsto \langle g, h \rangle,$ satisfying the
		following conditions:
		\begin{enumerate}
			\item $\langle g, hh' \rangle  = \langle g, h\rangle \langle ^hg, {^hh'}\rangle$,
			\item $\langle gg', h \rangle  = \langle ^gg', {^gh} \rangle \langle g, h\rangle$,
			\item $\langle (g \star g'), {^{g'}h} \rangle \langle ^hg, \langle g', h \rangle \rangle ^{-1} \langle ^gg', {\langle g, h \rangle} ^{-1}\rangle^{-1} = 1$,
			\item $\langle ^{h'}g, (h\star h')\rangle  ({^hh'}\star {\langle g, h\rangle})   ({^gh}\star  {\langle g, h'\rangle^{-1}}) = 1$,  
		\end{enumerate}
		where $g, g' \in G, h, h'\in H$, ${^gh}$ and ${^hg}$ denote group actions of $g$ on $h$ and $h$ on $g$, respectively and ${^gg'} = gg'g^{-1}, {^hh'} = hh'h^{-1}$.
	\end{definition}
	
	\begin{definition}\cite{DASS}
		Let $G$ and $H$ be multiplicative Lie algebras acting on each other. The actions are said to be compatible if:
		\begin{enumerate}
			\item $^{(^g h)}g' =\ ^{ghg^{-1}}g' \quad \text{and} \quad ^{(^hg)}h' =\ ^{hgh^{-1}}h'$;
			\item $\langle \langle h, g \rangle^{-1}, h' \rangle = \langle g, h \rangle \star h' \quad \text{and} \quad \langle \langle g, h \rangle^{-1}, g' \rangle = \langle h, g \rangle \star g'$;
			\item $^{\langle g, h \rangle \langle h, g \rangle} h' = h' \quad \text{and} \quad ^{\langle g, h \rangle \langle h, g \rangle} g' = g'$;
			\item $^{g} \langle h, g' \rangle = \langle ^gh, {^g g'} \rangle \quad \text{and} \quad ^{h} \langle g, h' \rangle = \langle {^hg}, {^hh'} \rangle$;
			\item $\langle g ^{h}g^{-1}, h' \rangle = (^ghh^{-1}) \star h' \quad \text{and} \quad \langle h^{g}h^{-1}, g' \rangle = (^hg g^{-1}) \star g'$,
		\end{enumerate}
		for all $g, g' \in G$ and $h, h' \in H$.
	\end{definition}
	\begin{lemma}\cite{DASS}
		Let \( G \) and \( H \) be multiplicative Lie algebras acting on each other compatibly. Then for all \( g \in G \) and \( h, h' \in H \),
		\[
		[\langle g,h\rangle, h'] = \langle g^{h}g^{-1}, h' \rangle = (^ghh^{-1}) \star h'.
		\]
	\end{lemma}
	
	
	\section{Relative Lie central extension and Frattini Subalgebra in multiplicative Lie algebras}
	
	In this section, we define the notion of relative Lie central extension and Frattini subalgebra for multiplicative Lie algebras.
	
	\begin{definition}\label{relative_def}
		A relative Lie central extension of the pair of multiplicative Lie algebra $(H, G)$ is a multiplicative Lie algebra homomorphism $\tau : L \to G$ where $G$ and $L$ compatibly acting on each other, and satisfying the following conditions
		\begin{enumerate}
			\item $\tau(L) = H$ and $G$ acts trivially on $\ker \tau$ 
			\item $\tau(^gl) = g \tau(l) g^{-1}$ and $\tau(\langle g, l \rangle) = g \star \tau(l)$
			\item $^{\tau(l)}l' = l l' l^{-1}$ and $\langle \tau(l), l' \rangle = l \star l'$
			\item ${^{l}g} = {^{\tau{(l)}}g}$ and $\langle l, g \rangle = \tau(l) \star g$
		\end{enumerate}
		for all $l, l' \in L$ and $g \in G$.
	\end{definition}
	\textbf{Note:} By a pair $(H, G)$ of multiplicative Lie algebra, we mean a multiplicative Lie algebra $G$ and an ideal $H$ of $G$.
	
	\begin{definition}\label{pair_homomorphism}
		
		Let $\tau_i : L_i \to G_i$, $(i = 1,2)$ be two relative Lie central extensions of the multiplicative Lie algebra pairs $(H_1, G_1)$ and $(H_2, G_2)$. If there exist multiplicative Lie algebra homomorphisms $\theta : G_1 \to G_2$ with $\theta(H_1) = H_2$ and $\beta : L_1 \to L_2$ with $\beta(^M\{L_1, G_1\}) = {^M\{L_2, G_2\}}$ such that $\theta \tau_1(g_1) = \tau_2 \beta(g_1)$ for all $g_1 \in G_1$. Then the pair $(\theta, \beta) : \tau_1 \to \tau_2$ is called a morphism from $\tau_1$ to $\tau_2$.
		
	\end{definition}
	
	
	\begin{definition}\label{com_center_def}
		Lie commutator and Lie center for a pair of multiplicative Lie algebras $(H, G)$ are defined as, respectively
		\begin{enumerate}
			\item $^M[G, H] = [G, H](G \star H) = \langle [g, h] (g'\star h') : g, g' \in G,~ h, h' \in H \rangle$,
			\item $\mathcal{Z}(G, H) = \{h \in H : [g, h] = 1,~ g \star h = 1 ~ \text{for all} ~ g \in G\}$.
		\end{enumerate}
	\end{definition}
	
	
	\begin{definition}\label{rel_com_def}
		Let multiplicative Lie algebra $L$ and $G$ acting on each other compatibly. Then $G$-Lie commutator and $G$-Lie center of $L$ are defined as, respectively
		\begin{enumerate}
			\item $^{M}\{G, L\} = \langle ^gll^{-1} \langle g', l'\rangle : l, l' \in L,~ g, g' \in G \rangle$,
			\item $\bar{\mathcal{Z}}(G, L) = \{ l \in L : {^gl} = l,~ \langle g, l\rangle = 1 ~ \text{for all} ~ g \in G\}$.
		\end{enumerate}
	\end{definition}
	A pair $(H, G)$ of multiplicative Lie algebras is called Lie perfect pair if $^M\{G, H\} = H$. In particular, $^M[G, H] = H$.
	\begin{remark}
		\begin{enumerate}
			\item Using identities $(4)$ and $(5)$ of the Definition \ref{relative_def}, $\mathrm{Ker}(\tau) \subseteq \mathcal{Z}(L)$ and also seen that $\mathrm{Ker}(\tau) \subseteq \bar{\mathcal{Z}}(G, L)$. Let $H = G$, then the relative Lie central extension $\tau : L \to G$ gives the following Lie central extension of multiplicative Lie algebras
			\[\tilde{E} \equiv \xymatrix{1 \ar[r] & \mathrm{Ker} (\tau) \ar[r]^{} & L \ar[r]^{\tau} & G \ar[r] &1.}\]
			\item The $G$-Lie commutator and the $G$-Lie center of $L$ are ideals in $L$ \cite{DK_crossed}.
		\end{enumerate}
		\end{remark}
	Here, we introduce the concept of the Frattini subalgebra of a multiplicative Lie algebra and discuss its various properties.
	\begin{definition}\label{frattini_def}
		Let $G$ be a multiplicative Lie algebra. If $G$ has maximal subalgebra, then the intersection of all maximal subalgebras of $G$ is defined $\mathcal{M}$-Frattini subalgbera and denoted by $\tilde{{\Phi}}(G)$. Otherwise we defined $\tilde{\Phi}(G) = G$, if $G$ has no maximal subalgebra.
	\end{definition}
	
	\begin{examples}\label{examples_frattini}
		\begin{enumerate}
			\item 	Let $G = V_4=\langle a,b :a^2=b^2=1, ab=ba \rangle$ and $\star$ on $V_4$ is given by $a\star b=a$. Then maximal subalgebras of $V_4$ are $\{1, a\}, \{1, b\}$ and $\{1, ab\}$. Thus $\tilde{{\Phi}}(G) = 1$. 
			\item 	Let $G= D_4 = \langle a,b: b^4 = a^2=1 , ab= b^{-1} a \rangle$ and $\star$ on $D_4$ is given by $a \star b = b.$ Then maximal subalgebras of $D_4$ are $\{1, b, b^2, b^3\}, \{1, b^2, a, b^2a\}$ and $\{1, b^2, ba, b^3a\}$.  Thus $\tilde{\Phi}(G) = \{1, b^2\}$.
			\item Let $G = S_3$ and $\star$ is defined as $a \star b = [a, b]$ for all $a, b \in S_3$. Then maximal subalgebras of $S_3$ are $A_3$ and isomorphic to $\{1, (12)\}$. Thus $\tilde{{\Phi}}(S_3)= 1$.
		\end{enumerate}
	\end{examples}
	
	\begin{definition}\label{non-generator}
		Let $G$ be a multiplicative Lie algebra. An element $g \in G$ is called non-generator of G, if $G = \langle X \rangle$, i.e. $G$ is generated by some subset $X$ of $G$, then $G = \langle X - \{g\} \rangle
		$.
	\end{definition}
	
	\begin{proposition}\label{Frattini_Non}
		Let $G$ be a multiplicative Lie algebra. Then the $\mathcal{M}$-Frattini subalgebra $\tilde{\Phi}(G)$ is equal to the set of all non-generators of $G$.
		\begin{proof}
			Let $Z$ be the set of all non-generators of $G$. Suppose $g \notin Z$, then there exists a $X \subseteq G$ such that $ \langle X \cup \{g\} \rangle = G $ and $\langle X \rangle \neq G$. Let $K = \langle X \rangle$ and $\mathcal{F} = \{ L \subsetneq G : K \subseteq L \}$. Clearly, $\mathcal{F}$ is a nonempty partially ordered set with respect to inclusion. Let $\{J_\alpha : \alpha \in I\}$ be a chain of subalgebras of $G$ in $\mathcal{F}$. Then $\underset{\alpha \in I}{\cup} J_\alpha$ is a subalgebra of $G$, also $g \notin \underset{\alpha \in I}{\cup} J_\alpha$, for $g \notin J_\alpha$ for all $\alpha \in I$. Hence $\underset{\alpha \in I}{\cup} J_\alpha \in \mathcal{F}$. Therefore, every chain in $\mathcal{F}$ has an upper bound, by Zorn's lemma, $\mathcal{F}$ has a maximal element $M$. Then $g \notin M$, if $N$ is a subalgebra of $G$ such that $M \subsetneq M$, then $N \notin \mathcal{F}$, this implies $g \in N$ and $G = N$. Thus $M$ is maximal subalgebra of $G$. Hence $x \notin \tilde{\Phi}(G)$. Therefore $\tilde{\Phi}(G) \subseteq Z$.
			
			Conversely, let $g \notin \tilde{{\Phi}}(G)$, so there exists a maximal subalgebra $M$ of $G$ such that $g \notin M$. Then $\langle M \cup \{x\} \rangle = G$ and $\langle M \rangle = M \neq G$. Thus $x \notin X$ and $X \subseteq \tilde{{\Phi}}(G)$.
		\end{proof}
	\end{proposition}
	
	\begin{definition} \label{normalizer_MLA}
		Let $H$ be a subalgebra of multiplicative Lie algebra $G$. Then the multiplicative normalizer of $H$ is defined as 
		\[\mathcal{N}_G(H) = \{g \in G : g H g^{-1} = H ~ \text{and}~ g \star H \subseteq H\}.\]
		
	\end{definition}
	
	\begin{proposition}\label{normalizer_subalgebra}
		Let $H$ be a subalgebra of multiplicative Lie algebra $G$. Then $\mathcal{N}_G(H)$ is also a subalgebra of $G$.
		\begin{proof}
			Let $g_1, g_2 \in \mathcal{N}_G(H)$. Then $^{g_i}h$ and $g_i \star h \in H$ for all $h \in H$ and $i = 1, 2$. Since $g_1 g_2 \star h = {^{g_1}(g_2 \star h)} (g_1 \star h)$ and $g_1^{-1} \star h = {^{g_1^{-1}}(g_1 \star h) },$ this implies $g_1 g_2 \in \mathcal{N}_G(H)$ and $g_1^{-1} \in \mathcal{N}_G(H)$. Thus $\mathcal{N}_G(H)$ is a subgroup of $G$. Now consider
			\begin{align*}
				((g_1 \star g_2) \star h)^{-1} & = \big((g_1 \star g_2) \star {^{^{g_2}g_2^{-1}}h}\big)^{-1}\\
				& = \big((g_2 \star {^{g_2^{-1}}}h) \star {^{^{g_2^{-1}}h}g_1}\big) \big((^{g_2^{-1}}h \star g_1) \star {^{g_1} g_2}\big)
			\end{align*}
			First consider
			\begin{align*}
				\big((g_2 \star {^{g_2^{-1}}}h) \star {^{^{g_2^{-1}}h}g_1}\big) & = \big(h' \star (^{g_2^{-1}}h) g_1 (^{g_2^{-1}}h)^{-1}\big), ~ \text{where} ~ h' = (g_2 \star {^{g_2^{-1}}}h)\\
				& = \big(h' \star {^{g_2^{-1}}h}\big) {^{^{g_2^{-1}}h} \big(( h' \star g_1) {^{g_1}(h' \star (^{g_2^{-1}}h)^{-1})}\big)}.
			\end{align*}
			Thus $\big((g_2 \star {^{g_2^{-1}}}h) \star {^{^{g_2^{-1}}h}g_1}\big) \in H$ and similarly $\big((^{g_2^{-1}}h \star g_1) \star {^{g_1} g_2}\big) \in H$, this implies $g_1 \star g_2 \in \mathcal{N}_G(H)$. Hence $\mathcal{N}_G(H)$ is a subalgebra of $G$.
		\end{proof}
	\end{proposition}
	
	
	Let $H$ be a subalgebra of multiplicative Lie algebra $G$. Define series of subalgebras $ \mathcal{N}_G^n(H) $, where $n \in \mathbb{N} ~ \cup \{0\} $ of $G$, inductively as $\mathcal{N}_G^0(H) = H$, assuming that $\mathcal{N}_G^n(H)$ has been already defined, then define $\mathcal{N}_G^{n+1}(H) = \mathcal{N}_G(\mathcal{N}_G^n(H))$. Thus we get an ascending chain
	\[ H \trianglelefteq   \mathcal{N}_G^1(H)    \trianglelefteq      \mathcal{N}_G^2(H)    \trianglelefteq \cdots \trianglelefteq  \mathcal{N}_G^n(H) \trianglelefteq \cdots . \]
	
	
	\begin{proposition}\label{Normalizer= G}
		Let $G$ be a nilpotent multiplicative Lie algebra of class $n$. Then $ \mathcal{N}_G^n(H) = G$ for all subalgebras $H$ of $G$.
		\begin{proof}
			Assume $\mathcal{Z}_n(G) = G$. We will show by induction on $r$ that $\mathcal{Z}_r(G) \subseteq \mathcal{N}_G^r(H)$. Let $r =0 $. Then $\mathcal{Z}_0(G) = \{1\} \subseteq H = \mathcal{N}_G^0(H)$. Suppose $r = k $ and $\mathcal{Z}_k(G) \subseteq \mathcal{N}_G^k(H)$. Now, let $r = k+1$ and $g \in \mathcal{Z}_{k+1}(G)$, then $g \mathcal{Z}_k(G) \in \mathcal{Z}\Big(\dfrac{G}{\mathcal{Z}_k(G)} \Big)$. Let $x \mathcal{Z}_k(G) \in \dfrac{G}{\mathcal{Z}_k(G)}$. Then $g x g^{-1} x^{-1}$ and $g \star x \in \mathcal{Z}_k(G)$ for all $x \in G$. If $x \in \mathcal{N}_G^k(H)$, then $g x g^{-1}$ and $g \star x \in \mathcal{N}_G^k(H)$, this implies $g \in \mathcal{N}_G^{k+1}(H)$. Thus $\mathcal{Z}_{r}(G) \subseteq \mathcal{N}_G^{r}(H)$ for all $r$. Hence $G = \mathcal{Z}_n(G) \subseteq \mathcal{N}_G^n(H)$. Therefore $\mathcal{N}_G^n(H) = G$.         	\end{proof}
	\end{proposition}
	
	
	\begin{definition}\label{normalizer_condition}
		A multiplicative Lie algebra is said to satisfy the multiplicative normalizer condition if every proper subalgebra $H$ of $G$ is properly contained in its multiplicative normalizer, i.e. $H \subsetneq \mathcal{N}_G(H)$.
	\end{definition}
	
	
	\begin{corollary}\label{corollary_norm_condition}
		Every nilpotent multiplicative Lie algebra $G$ satisfy the multiplicative normalizer condition.
		\begin{proof}
			Let $G$ be a nilpotent multiplicative Lie algebra of class $n$ and $H$ be a proper subalgebra of $G$. Then $\mathcal{N}_G^n(H) = G$ (by proposition \ref{Normalizer= G}). Suppose $H = \mathcal{N}_G(H)$, this implies that $\mathcal{N}_G^n(H) = H \neq G$. Which is a contradiction.
		\end{proof}
	\end{corollary}
	
	\begin{remark}\label{nilpotent_remark}
		In general, a maximal subalgebra is not necessarily an ideal (see Example \ref{examples_frattini} (3)). However, in the following corollary, we demonstrate that this holds for nilpotent multiplicative Lie algebras.
	\end{remark}
	
	\begin{corollary}\label{maximal_are_ideals}
		Let $G$ be a nilpotent multiplicative Lie algebra. Then maximal subalgebras are ideals in $G$.
		\begin{proof}
			Let $M$ be a maximal subalgebra of $G$. Then $M \subsetneq \mathcal{N}_G(M)$ (by corollary \ref{corollary_norm_condition}), this implies $\mathcal{N}_G(M) = G$. Hence $M$ is an ideal of $G$. 
		\end{proof}
	\end{corollary}

	\begin{proposition}\label{derived_implies_equal}
		If $G$ is a nilpotent multiplicative Lie algebra of class $n$ and $H$ is a subalgebra of $G$ such that $G = {^M[G, G]}H$. Then $H = G$.
	\end{proposition}
	\begin{proof}
		Suppose $H \neq G$. Then by Proposition \ref{Normalizer= G}, $\mathcal{N}_G^0(H) = H,~ \mathcal{N}_G^{n-1}(H) \neq G$, but $\mathcal{N}_G^n(H) = G$, this implies that $\mathcal{N}_G^{n-1}(H)$ is a proper ideal of $G$ and $G/\mathcal{N}_G^{n-1}(H)$ is an abelian multiplicative Lie algebra with trivial Lie product for $G/\mathcal{Z}_{n-1}(G)$ is an abelian multiplicative Lie algebra with trivial Lie product and $\mathcal{Z}_{n-1}(G) \subseteq \mathcal{N}_G^{n-1}(H)$. Thus $^M[G, G] \subseteq \mathcal{N}_G^{n-1} (H)$. This implies $^M[G, G]H \subseteq {^M[G, G]}\mathcal{N}_G^{n-1} (H) = \mathcal{N}_G^{n-1}(H) \neq G$. Which is a contradiction. Hence $H = G$.
	\end{proof}
	 The Lie product $\star$ on $G$ is called \emph{proper} if there exist elements $r, s, g, g' \in G$ such that $r \star s \neq 1$ and $g \star g' \neq [g, g']$. A group $G$ is said to be \emph{Lie simple} \cite{RLS} if it admits no proper Lie products.
	\begin{remark}\label{no proper subalgebra}
		If multiplicative Lie algebra $G$ has no proper subalgebra then $G$ is cyclic group of prime order. If $G= \langle g \rangle$ is infinite cyclic group then $G$ is Lie simple, consider subgroup $H = \langle g^2 \rangle$ of $G$, then $H$ is a proper subalgebra of $G$, which is a contradiction. If $G$ is non-cyclic group of finite order, then there exist a proper subalgebra $H = \langle h \rangle $ for some $h \in G$, which is again a contradiction. Similarly, we can check the other cases.
	\end{remark}

	\begin{corollary}\label{derived_subset}
		Let $G$ be a nilpotent multiplicative Lie algebra. Then $$^M[G, G] \subseteq \tilde{{\Phi}}(G).$$
	\end{corollary}
	\begin{proof}
		
		Let $M$ be a maximal subalgebra of $G$. Since $G$ is nilpotent, Corollary \ref{maximal_are_ideals} implies that $M$ is a maximal ideal of $G$. Thus, Remark \ref{no proper subalgebra}, gives $G/M$ is an abelian multiplicative Lie algebra with trivial Lie product. Hence $^M[G, G] \subseteq \tilde{{\Phi}}(G)$.
	\end{proof}
	
	\begin{remark}
		If $H$ is a subalgebra of $G$ and $K$ is an ideal of G. Then $HK$ generated by the set $\{hk : h \in H~ \text{and} ~ k \in K\}$ is a subalgebra of $G$.
	\end{remark}

	\begin{proposition}\label{intersection_subset_of_frattini}
		Let $G$ be a multiplicative Lie algebra. Then $$\mathcal{Z}(G) \cap {^M[G, G]}  \subseteq \tilde{\Phi}(G).$$
		\begin{proof}
			Let $K = \mathcal{Z}(G) \cap {^M[G, G]} $ and suppose that $K \nsubseteq \tilde{{\Phi}}(G)$. Then there exist a maximal subalgebra $M$ of $G$ such that $K \nsubseteq M.$ But $M \leq MK \leq G$, implies $MK = G$. Let $m \in M$ and $g \in G$ then $g = m'k$ for some $m' \in M$ and $k \in K$. Clearly, $g m g^{-1} \in M$ and $m \star g \in M$ for $m \star m'k = (m \star m'){^{m'}(m \star k)} = m \star m'$, $(\because ~ k \in \mathcal{Z}(G))$. Thus $M$ is an ideal of  $G$. Hence $G/M$ is an abelian multiplicative Lie algebra with trivial Lie product. Therefore $^M[G, G] \leq M, ~ K \leq M$. Which is a contradiction. 
		\end{proof}
	\end{proposition}

	\begin{lemma}\label{lemma_identities}
		Let $G$ be a finite multiplicative Lie algebra and $K$ be an ideal of $G$. Then following are true:
		\begin{enumerate}
			\item $K  \leq \tilde{{\Phi}}(G)$ if and only if there is no proper subalgebra of $G$ such that $HK = G$.
			\item If $H \leq G$ and $K \leq \tilde{{\Phi}}(H)$ then $K \leq \tilde{{\Phi}}(G)$.
			\item If $K$ is nilpotent then $\tilde{{\Phi}}(K) \leq \tilde{{\Phi}}(G)$. 
		\end{enumerate}
		\begin{proof}
			\textit{(1)~} 
			Let $H < G$. Then there exist a maximal subalgebra $M$ of $G$ such that $HK \leq M < G$. Thus $HK \leq M < G$. Hence there is no proper subalgebra of $G$ such that $HK = G$. Conversely, suppose that $K \nleq \tilde{{\Phi}}(G)$. Then there exist a maximal subalgebra $M$ of $G$ such that $K \nleq M$. Thus $M < MK \leq G$, this implies $MK = G$. Which is a contradiction.\\
			\textit{(2)~} Suppose $K \nleq \tilde{{\Phi}}(G)$. Then there exist a proper subalgebra $J$ of $G$ such that $JK = G$. Since $K \leq H \leq G = JK$, this imples $H = H \cap JK = (H \cap J) K$. Thus $H = H \cap J$ for $K \leq \tilde{{\Phi}}(H)$. Hence $K \leq H \leq J$. Therefore $G = JK = J$. Which is a contradiction.\\
			\textit{(3)~} By Corollary \ref{maximal_are_ideals}, $\tilde{{\Phi}}(K)  \trianglelefteq K$, this implies $\tilde{{\Phi}}(K)  \trianglelefteq G$. Now replacing $K$ by $\tilde{{\Phi}}(K)$ and $H$ by $K$ in the last part, gives $\tilde{{\Phi}}(K) \leq \tilde{{\Phi}}(G)$.
		\end{proof}
	\end{lemma}


	\section{Isoclinism of relative Lie central extensions}
	
	In this section, we define the concept of isoclinism for relative Lie central extensions of multiplicative Lie algebras and prove some results related to the isoclinism.
	
	
	\begin{definition}\label{isoclinism_relative}
		
		Let $\tau_i : L_i \to G_i$, $i = 1,2$ be two relative multiplicative central extensions of the multiplicative Lie pairs $(G_i, L_i)$, $(i = 1, 2)$, respectively. Then the relative Lie central extension $\tau_1$ and $\tau_2$ are said to be isoclinic if there exist multiplicative Lie algebra isomorphisms $\theta : G_1  \to G_2$ with $\theta(H_1) = H_2$ and $\beta : {^M\{G_1, L_1\}} \to {^M\{G_2, L_2\}}$ such that
		\[\beta (^{g_1}l_1l_1^{-1}) = {^{g_2}l_2}l_2^{-1} ~ \text{and} ~ \beta(\langle g_1', l_1' \rangle) = \langle g_2', l_2' \rangle\]
		for all $g_1, g_1' \in G_1$ and $l_1, l_1' \in L_1$, where $g_2, g_2' \in G_2$ and $l_2, l_2' \in L_2$ such that	\[\theta \tau_1(l_1) = \tau_2(l_2),~ \theta \tau_1(l_1') = \tau_2(l_2'),~ \theta(g_1) = g_2 ~ \text{and}~ \theta(g_1') = g_2'.\]
		We called the pair $(\theta, \beta)$ to be isoclinism from $\tau_1$ to $\tau_2$ and denoted as $(\theta, \beta) : \tau_1 \sim^{mlr} \tau_2$.
		
	\end{definition}
	
	\begin{definition}\label{isoclinism_2_def}
		A morphism $(\theta, \beta) : \tau_1 \to \tau_2$ is called isoclinic morphism if the pair $(\theta, \beta|_{^M\{G_1, L_1\}})$ is an isoclinism from $\tau_1$ to $\tau_2$. Also, $(\theta, \beta)$ is isoclinic monomorphism and epimorphism if $\beta$ is monomorphism and epimorphism, respectively.
	\end{definition}
	
	Now, we will give some examples of isoclinism of relative Lie central extensions.
\begin{examples}\label{examples_iso}
	\begin{enumerate}
		\item Let $K$ be an abelian multiplicative Lie algebra with trivial Lie bracket and $L$ be a multiplicative Lie algebra. Also, let $\phi_{L} : L \times K \to L$, $i_{L} : L \to L \times K$ be projection and canonical homomorphisms, respectively. \\
		Suppose $\tau : L \to G$ is a relative Lie central extension of the pair $(H,G)$. Then $\tau \phi_{L} : L \times K \to G$, $\tau \phi_{L} (l, k) = \tau(l)$ is a multiplicative Lie algebra homomorphism. Suppose action of $G$ on $L \times K$ is defined as
		\[^g(l, k) = (^gl, k); \hspace{1 cm} \langle g, (l, k) \rangle = (\langle g, l \rangle, 1 ),\]
		and action of $L \times K$ on $G$ is defined as 
		\[{^{(l,k)}g} = {^lg}; \hspace{1 cm} \langle (l, k), g \rangle = \langle l, g \rangle.\]
		Then, $G$ and $L \times K$ acting on each other combatibly and morphism $\tau \phi_{L}$ is relative Lie central extension. 
		\begin{align*}
			\hspace{1 cm}	^M\{G, L \times K\} & = \langle ^g(l,k) (l,k)^{-1} \langle g', (l',k') \rangle : g, g' \in G,~ (l,k), (l',k') \in L \times K \rangle \\
			& = \langle (^gl l^{-1}\langle g', l'\rangle, 1 ) : g, g' \in G, ~ l, l' \in L \rangle \cong {^M\{G, L\}}.
		\end{align*}
		Clearly, the morphisms $(\phi_{L}, I_G) : \tau \phi_{L} \to \tau$ and $(i_{L}, I_{L}) : \tau \to \tau \phi_{L}$ are isoclinic epimorphism and isoclinic monomorphism, respectively.
		\item Let $M$ be a subalgebra of $L$ such that $L = MH$ and $\tau(M) = H$, where $\tau : L \to G$ is a relative Lie central extension of the pair $(H, G)$. Then $\tau|_M : M \to G$ defined by $\tau|_M (m) = \tau(m)$ for all $m \in M$ is a relative Lie central extension of $(H, G)$. Also, the morphism $(i, I_G) : \tau|_G \to \tau$ is isoclinic monomorphism, where $i$ is inclusion map.
		
		\item Let $\tau : L \to G$ be a relative Lie central extension of the pair $(H, G)$ and $K \subseteq \mathrm{Ker}(\tau)$. Then $\bar{\tau} : L/ K \to G$ defined by $\bar{\tau}(m K) = \tau(m)$ is a morphism. Suppose the action of $G$ on $L/K$ is defined as
		\[(g, mK) \mapsto {^g({m K})} = {^gm}K; \hspace{1 cm} (g, mK) \mapsto \langle g, mK \rangle = \langle g, m \rangle K,\]
		and the action of $L/K$ on $G$ is defined as 
		\[(mK, g) \mapsto {^{mK}g} = {^mg}; \hspace{1 cm} (mK, g) \mapsto \langle mK, g \rangle = \langle m, g \rangle.\]
		Let $(mK, g) = (m'K, g')$. Then $m = m'k$ for some $k \in K$, compute 
		\[^g mK = {^gm} K = {^g}m'k K = {^gm'} {^gk}K = {^{g'}m'}K;~ (\because~ K \subseteq \mathrm{Ker}(\tau)),\]
		\[\hspace{1cm} \langle g, m \rangle K = \langle g, m'k \rangle K = \langle g', m' \rangle K \langle {^{m'}g'}, {^{m'}k} \rangle K = \langle g', m' \rangle K,\]
		\[^{mK} g = {^mg} = {^{\tau(m)}}g = {^{\tau(m'k)}}g = {^{\tau(m')}(^{\tau(k)}g)} = {^{\tau(m')}}g'= {^{m'}}g',\]
		\[\langle mK, g \rangle = \langle m'k, g \rangle = {^{m'}\langle k, g \rangle} \langle m', g \rangle = {^{m'}(\tau(k) \star g)} \langle m', g \rangle = \langle m', g \rangle.\]
		Thus these actions are well-defined. It can be easily check that $G$ and $L/K$ acting on each other compatibly and $\bar{\tau}$ is relative Lie central extension. Now, let $\gamma : L \to L/I$ be a natural epimorphism. Then $(\gamma, I_G) : \tau \to \bar{\tau}$ is an isoclinic epimorphism.
	\end{enumerate}
\end{examples}


\begin{lemma}\label{imp_lemma}
	Suppose $\tau_1 : L_1 \to G_1$ and $\tau_2 : L_2 \to G_2$ are two relative Lie central extensions of the pairs $(H_1, G_1)$ and $(H_2, G_2)$, respectively. If $(\theta, \beta) : \tau_1 \to \tau_2$ is an relative isoclinism, then the following are true
	\begin{enumerate}
		\item $\theta \tau_1(x) = \tau_2 \beta(x)$ for all $x \in {^M\{G_1, L_1\}}$
		\item $\beta(\ker \tau_1 \cap {^M\{G_1, L_1\}}) = \ker \tau_2 \cap {^M\{G_2, L_2\}}$
		\item $\beta(^{g_1}l_1 l_1^{-1} \langle g_1', l_1' \rangle) = {^{g_2} \beta(l_1)} \beta(l_1)^{-1} \langle g_2', \beta(l_1') \rangle$, for all $l_1, l_1' \in {^M\{G_1, L_1\}}$ with $\theta(g_1) = g_2$ and $\theta(g_1') = g_2'$.
	\end{enumerate}
\end{lemma}
\begin{proof}
	\textit{(1)} Let $x = {^{g_1}l_1 l_1^{-1}} \langle g_1', l_1' \rangle \in {^M\{G_1, L_1\}}$ for $l_1, l_1' \in L_1$ and $g_1, g_1' \in G_1$. Then 
	\begin{align*}
		\theta \tau_1(  {^{g_1}l_1 l_1^{-1}} \langle g_1', l_1' \rangle ) & = \theta (\tau_1(^{g_1}l_1) \tau_1(l_1)^{-1}) \theta (\tau_1 \langle g_1', l_1' \rangle) \\
		& = \theta (g_1 \tau_1(l_1) g_1^{-1} \tau_1(l_1)^{-1}) (\theta (g_1') \star \theta (\tau_1(l_1')))
	\end{align*}
	If $\theta(g_1) = g_2, \theta(g_1') = g_2'$ and $\theta(\tau_1(l_1)) = \tau_2(l_2)$, $\theta(\tau_1(l_1')) = \tau_2(l_2')$ for $g_2, g_2' \in G_2$ and $l_2, l_2' \in L_2$. Then
	\begin{align*}
		\theta \tau_1(  {^{g_1}l_1 l_1^{-1}} \langle g_1', l_1' \rangle ) & = [g_2, \tau_2(l_2)] (g_2' \star \tau_2(l_2')) \\
		& = \tau_2(^{g_2}l_2 ) \tau_2(l_2)^{-1} \tau_2(\langle g_2', l_2' \rangle) \\
		& = \tau_2(^{g_2}l_2 l_2^{-1} \langle g_2', l_2' \rangle) \\
		& = \tau_2(\beta(^{g_1}l_1 l_1^{-1} \langle g_1', l_1' \rangle)) = \tau_2(\beta (x)).
	\end{align*}
	
	\textit{(2)} Let $x \in  \ker \tau_1 \cap {^M\{G_1, L_1\}}$, then $\tau_1(x) = 0$ and $x = {^{g_1}l_1}l_1^{-1} \langle g_1', l_1' \rangle$ for some $l_1, l_1' \in L_1$ and $g_1, g_1' \in G$. By Definition, \ref{isoclinism_relative} $\beta(^{g_1} l_1 l_1^{-1}) = {^{g_2}l_2} l_2^{-1}$ and $\beta(\langle g_1', l_1' \rangle) = \langle g_2', l_2' \rangle$ for $\theta(g_1) = g_2, \theta(g_1') = g_2', \theta \tau_1(l_1) = \tau_2(l_2)$ and $\theta \tau_1(l_1') = \tau_2(l_2')$, compute
	\[\tau_2(^{g_2}l_2 l_2^{-1}\langle g_2', l_2' \rangle) = \tau_2(^{g_2}l_2) \tau_2(l_2)^{-1} \tau_2( \langle g_2', l_2' \rangle ) = [g_2, \tau_2(l_2)](g_2' \star \tau_2(l_2'))= \theta \tau_1(x) = 0.\]
	This implies $\beta(x) \in \ker \tau_2 \cap {^M\{G_2, L_2\}}$. Similarly, we can show the converse.
	
	\textit{(3)} Let $l_1, l_1' \in {^M\{G_1, L_1\}}$, $\theta(g_1) = g_2, \theta(g_1') = g_2'$ and $\theta \tau_1(l_1) = \tau_2(\beta(l_1)), \theta \tau_1(l_1') = \tau_2 \beta(l_1')$. By Definition \ref{isoclinism_relative}, we have $\beta(^{g_1}l_1 l_1^{-1}) = {^{g_2}\beta(l_1)}\beta(l_1)^{-1}$, $\beta(\langle g_1', l_1' \rangle ) = \langle g_2', \beta(l_1') \rangle$.
\end{proof}


Let $\tau_1 : L_1 \to G_1$ and $\tau_2 : L_2 \to G_2$ be two relative Lie central extensions of the multiplicative Lie algebra pairs $(H_1, G_1)$ and $(H_2, G_2)$, respectively. Let $\theta : G_1\to G_2$ be a multiplicative Lie algebra. Define
\[\mathcal{L} = \{(l_1, l_2) : l_1 \in L_1, l_2 \in L_2 ~ \text{and}~ \theta \tau_1(l_1) = \tau_2(l_2) \}.\]
Then $\mathcal{L}$ is a subalgebra of $L_1 \times L_2$.


\begin{lemma}\label{lemma_ker}
	Let $\tau : \mathcal{L} \to G_1$, defined by $\tau(l_1, l_2) = \tau_1(l_1)$ be a multiplicative Lie algebra homomorphism. Then $G_1$ is acting on $\mathcal{L}$ with the actions defined as
	\[(g_1, (l_1, l_2)) \mapsto {^{g_1}(l_1, l_2)} = (^{g_1}l_1, {^{g_2}l_2}); \hspace{2 cm}  \langle g_1, (l_1,l_2) \rangle \mapsto (\langle g_1, l_1 \rangle , \langle g_2, l_2 \rangle ) \]
	and $\mathcal{L}$ is acting on $G_1$ with the actions defined as 
	\[	((l_1,l_2), g_1) \mapsto {^{(l_1, l_2)}}g_1 = {^{l_1}g_1}; \hspace{2 cm}  \langle (l_1,l_2), g_1 \rangle \mapsto \langle l_1, g_1 \rangle,\]
	where $\theta(g_1) = g_2$, and $\tau$ is a relative Lie central extension of the pair $(H_1, G_1)$.
\end{lemma}
\begin{proof}
	Let $g_1, g_1' \in G_1$ and $g_2, g_2' \in G_2$ with $\theta(g_1) = g_2$ and $\theta(g_1') = g_2'$, and $(l_1,l_2), (l_1',l_2') \in \mathcal{L}$. Then by Definition \ref{relative_def} (6), $^{l_1}g_1 = {^{\tau_1(l_1)}g_1}$ and $^{l_2}g_2 = {^{\tau_2(l_2)}g_2}$. Thus,  
	\[	\theta(^{l_1}g_1) = \theta (^{\tau_1(l_1)}g_1) = \theta \tau_1(l_1) \theta(g_1) \theta \tau_1(l_1)^{-1} = \tau_2(l_2) g_2 \tau_2(l_2)^{-1} = {^{\tau_2(l_2)}g_2} = {^{l_2}g_2}.\]
	Now, we will prove the conditions of angle bracket to prove that $G_1$ is acting on $\mathcal{L}$
	\begin{align*}
		\langle g_1, g_1', (l_1, l_2) \rangle &= (\langle g_1 g_1', l_1 \rangle, \langle g_2 g_2', l_2 \rangle) \\
		& = (\langle {^{g_1}g_1'}, {^{g_1}l_1} \rangle \langle g_1, l_1 \rangle, \langle {^{g_2}g_2'}, {^{g_2}l_2} \rangle \langle g_2, l_2 \rangle) \\
		& = ( \langle {^{g_1}g_1'}, {^{g_1}l_1} \rangle, \langle {^{g_2}g_2'}, {^{g_2}l_2} \rangle) (\langle g_1, l_1 \rangle, \langle g_2, l_2 \rangle)\\
		& = \langle {^{g_1}g_1'}, ({^{g_1}l_1}, {^{g_2}l_2}) \rangle \langle g_1, (l_1, l_2) \rangle =  \langle {^{g_1}g_1'}, {^{g_1}}(l_1, l_2) \rangle \langle g_1, (l_1, l_2) \rangle,
	\end{align*}
	\begin{align*}
		\langle g_1, (l_1, l_2) (l_1', l_2') \rangle & = \langle g_1, (l_1 l_1', l_2l_2') \rangle =( \langle g_1, l_1 l_1' \rangle , \langle g_1, l_2l_2' \rangle  ) \\
		& = (\langle g_1, l_1 \rangle \langle ^{l_1} g_1, {^{l_1}l_1'} \rangle, \langle g_2, l_2 \rangle \langle ^{l_2} g_2, {^{l_2}l_2'} \rangle) \\
		& = (\langle g_1, l_1 \rangle , \langle g_2, l_2 \rangle) (\langle ^{l_1} g_1, {^{l_1}l_1'} \rangle, \langle ^{l_2} g_2, {^{l_2}l_2'} \rangle) \\
		&  = \langle g_1, (l_1, l_2) \rangle \langle ^{l_1}g_1, (^{l_1}l_1', {^{l_2}l_2'}) \rangle = \langle g_1, (l_1, l_2) \rangle \langle ^{(l_1,l_2)}g_1, {^{(l_1,l_2)}}(l_1', l_2') \rangle,
	\end{align*}
	\begin{align*}
		\langle g_1 \star g_1', {^{g_1'}}& (l_1, l_2) \rangle \langle ^{(l_1,l_2)} g_1, \langle g_1', (l_1,l_2) \rangle \rangle ^{-1} \langle ^{g_1}g_1', \langle g_1, (l_1,l_2) \rangle ^{-1} \rangle ^{-1} \\
		& = \langle g_1 \star g_1', ({^{g_1'}}l_1, {^{g_2'}}l_2) \rangle \langle ^{l_1}g_1, (\langle g_1', l_1\rangle , \langle g_2', l_2 \rangle) \rangle^{-1} \langle ^{g_1}g_1', (\langle g_1, l_1 \rangle ^{-1} , \langle g_2, l_2 \rangle ^{-1}) \rangle ^{-1}\\
		& = (\langle g_1 \star g_1', {^{g_1'}}l_1 \rangle, \langle g_2 \star g_2', {^{g_2'}}l_2 \rangle) (\langle ^{l_1}g_1, \langle g_1', l_1 \rangle \rangle, \langle ^{l_2}g_2, \langle g_2', l_2 \rangle \rangle )^{-1} \\
		& ~~ \quad (\langle ^{g_1} g_1',  \langle g_1, l_1 \rangle ^{-1} \rangle , \langle ^{g_2} g_2',  \langle g_2, l_2 \rangle ^{-1} \rangle ) ^{-1} \\
		& = ( \langle g_1 \star g_1', {^{g_1'}}l_1 \rangle \langle ^{l_1}g_1, \langle g_1', l_1 \rangle \rangle^{-1} \langle ^{g_1} g_1',  \langle g_1, l_1 \rangle ^{-1} \rangle ^{-1}, \\
		& ~~~ \quad \langle g_2 \star g_2', {^{g_2'}}l_2 \rangle \langle ^{l_2}g_2, \langle g_2', l_2 \rangle \rangle ^{-1} \langle ^{g_2} g_2',  \langle g_2, l_2 \rangle ^{-1} \rangle ^{-1} ) = 1
	\end{align*}
	Similarly, we can show that
	\[\langle ^{(l_1', l_2')} g_1, (l_1, l_2) \star (l_1', l_2') \rangle (^{(l_1,l_2)}(l_1',l_2') \star \langle g_1, (l_1,l_2) \rangle ) (^{g_1}(l_1,l_2) \star \langle g_1, (l_1',l_2') \rangle ^{-1}) = 1. \]
	Clearly, $\mathcal{L}$ is acting on $G_1$ as this action induced by the action of $L_1$ on $G_1$. Also, the compatibilty of the actions can be easily check.
	Now, we will prove that $\tau$ is relative Lie central extension. Clearly, $\tau(\mathcal{L}) = H_1$ and
	\[\tau (^{g_1} (l_1, l_2)) = \tau (^{g_1}l_1, {^{g_2}}l_2)= \tau_1(^{g_1}l_1) = g_1 \tau_1(l_1) g_1^{-1} = g_1 \tau(l_1,l_2) g_1^{-1},\]
	\[\tau(\langle g_1, (l_1,l_2) \rangle) = \tau(\langle g_1, l_1\rangle , \langle g_2, l_2 \rangle) = \tau_1(\langle g_1, l_1 \rangle) = g_1 \star \tau_1(l_1) = g_1 \star \tau(l_1,l_2),\] 
	\[^{\tau(l_1,l_2)} (l_1', l_2') = (^{\tau_1(l_1)}l_1', {^{\tau_2(l_2)}}l_2') = (l_1,l_1'l_1^{-1}, l_2l_2'l_2^{-1}) = {^{(l_1,l_2)} (l_1',l_2')},\]
	\[\langle \tau(l_1,l_2) , (l_1',l_2') \rangle = (\langle \tau_1(l_1), l_1'\rangle, \langle \tau_2(l_2), l_2'\rangle ) = (l_1 \star l_1', l_2 \star l_2') = (l_1, l_2) \star (l_1', l_2'),\]
	\[^{\tau(l_1,l_2)} g_1 = {^{\tau_1(l_1)}} g_1 = {^{l_1}} g_1= {^{(l_1,l_2)} }g_1,\]
	\[ \tau(l_1, l_2)\star g_1 = \tau_1(l_1)\star g_1 = \langle l_1, g_1 \rangle = \langle (l_1, l_2), g_1 \rangle.\]
	Let $(l_1,l_2) \in \ker \tau.$ Then $l_1 \in \ker \tau_1$ and $l_2 \in \ker \tau_2$ for $\theta(\tau_1(l_1)) = \tau_2(l_2)$. Thus $^{g_1}(l_1,l_2) = ( ^{g_1}l_1, {^{g_2}}l_2) = (l_1, l_2)$ and $\langle g_1, (l_1, l_2 ) \rangle = ( \langle g_1, l_1 \rangle , \langle g_2 , l_2 \rangle) = 1$. Therefore, $\tau$ is a relative Lie central extension of the pair $(G_1, H_1)$.
\end{proof}

\begin{proposition}
	Let $\tau_1 : L_1 \to G_1$ and $\tau_2 : L_2 \to G_2$ be two relative Lie central extension of the pair $(H_1, G_1)$ and $(H_2, G_2)$, respectively. The isomorphism $\theta : G_1 \to G_2$ induces an isoclinism from $\tau_1$ to $\tau_2$ if and only if there exist isoclinic epimirohisms from $\tau$ onto $\tau_1$ and $\tau_2$.
\end{proposition}
\begin{proof}
	Let $(\theta_1, \beta_1)$ and $(\theta_2, \beta_2)$ be isoclinic epimorphisms, i.e. $(\theta_1, \beta_1|_{^M\{G_1, \mathcal{L}\}})$ and $(\theta_1, \beta_1|_{^M\{G_2, \mathcal{L}\}})$ are isoclinisms from $\tau$ to $\tau_1$ and $\tau_2$, respectively. Then $(\theta_2, \beta_2 \circ \beta_1^{-1})$ gives the isoclinism from $\tau_1$ to $\tau_2$.\\
	Conversely, let $\theta: G_1 \to G_2$ induces an isoclinism from $\tau_1$ to $\tau_2$. Then there exist an isomorphism $\beta : {^M\{G_1, L_1\}} \to {^M\{G_2, L_2\}}$. By Lemma \ref{imp_lemma} (1) gives 
	\[^M\{G_1, \mathcal{L}\} = \{(l_1, \beta(l_1)) : l_1 \in {^M\{G_1, L_1\}}\}.\]
	Let $\beta_i : \mathcal{L} \to L_i $, $(i = 1,2)$, defined as $\beta_1(l, \beta(l))= l$, $\beta_2(l, \beta(l))= \beta(l)$ and $\theta_1 = I_{G_1}$, $\theta_2 = \theta$. Clearly, $(\beta_1, \theta_1)$ and $(\beta_2, \theta_2)$ are isoclinic epimorphisms from $\tau$ to $\tau_1$ and $\tau_2$, respectively. 
\end{proof}
%
If $\mathcal{M} = \dfrac{\mathcal{L}}{^M\{G, \mathcal{L}\}}$, let $(l_1, l_2)^M\{G, \mathcal{L}\}, (l_1',l_2')^M\{G, \mathcal{L}\} \in \mathcal{M}$. Then by Lemma \ref{lemma_ker}, we have
\[[(l_1, l_2), (l_1',l_2')] = {^{\tau(l_1,l_2)}} (l_1', l_2')(l_1, l_2)^{-1} \in {^M\{G, \mathcal{L}\}},\]
\[(l_1, l_2) \star (l_1', l_2') = \langle \tau(l_1, l_2), (l_1', l_2') \rangle \in {^M\{G, \mathcal{L}\}}.\]
Thus $\mathcal{M}$ is abelian multiplicative Lie algebra with trivial Lie bracket. Suppose $\delta : \mathcal{L} \to G$ is a relative Lie central extension, then $\delta \phi_{\mathcal{L}} : \mathcal{L} \times \mathcal{M} \to G$ is also relative Lie central extension (see Example \ref{examples_iso} (1)). 
\begin{proposition}
	Let $(\theta, \beta)$ be an isoclinism from $\tau_1$ to $\tau_2$. Then there exist isoclinic monomorphisms from $\tau$ to $\tau_1 \phi_{L_1}$ and $\tau_2 \phi_{L_2}$.
\end{proposition}
\begin{proof}
	Let maps $\bar{\beta_i} : \mathcal{L} \to L_i \times M$ defined by $\bar{\beta_i}(n_1, n_2) = (\beta_i(n_1,n_2), (n_1, n_2)^M\{G_1, \mathcal{L}\})$, $(i=1,2)$, where $\beta_i$ are projection maps. We can easily verify that $\beta_i$ are multiplicative Lie algebra homomorphisms. If $(n_1, n_2) \in \ker \beta_1$, $n_1 = 1$ and $n_2 = 1$ for $n_2 = \beta(n_1)$, also if $(n_1, n_2) \in \ker \beta_2$, $n_2 = 1$ and $\beta(n_1) = n_2 = 1$ this implies $n_1= 1$ for $n_1 \in {^M\{G_1, L_1 \}}$. Thus $\ker \bar{\beta_i}= 1$, i.e. $\bar{\beta_i}$ are monomorphisms. If $\beta_i|_{^M\{G_i, \mathcal{L}\}}$, $(i= 1,2)$, then $\bar{\beta_1}(n_1, \beta(n_1)) = (n_1, {^M\{G_1, \mathcal{L}\}})$ and $\bar{\beta_2}(n_1, \beta(n_1)) = (\beta(n_1), {^M\{G_1, \mathcal{L}\}})$ are isomorphisms. Hence $(\bar{\beta_1}, I_{G_1}) : \bar{\tau} \to \tau_1 \phi_{L_1}$ and $(\bar{\beta_2}, \theta) : \bar{\tau} \to \tau_2 \phi_{L_2}$.
\end{proof}
\begin{proposition}
	Let $\theta : G_1 \to G_2$ be an isomorphism. Then $\theta$ induces an isoclinism from $\tau_1$ to $\tau_2$ if and only if there exist isoclinic monomorphisms from $\tau_1$ and $\tau_2$ to the relative Lie central extension $\overline{\tau_1 \phi_{L_1}} : (L_1 \times M)/N \to G_1$, where $N$ is an ideal of $L_1 \times M$ and $T \subseteq \ker \tau_1 \phi_{L_1}$.
\end{proposition}
\begin{proof}
	Suppose $\theta$ induces an isoclinism $(\theta, \beta)$ from $\tau_1$ to $\tau_2$. Let
	\[N = \{(l, (l, 1) ^M\{\mathcal{L}, G_1\}) : l \in \ker \tau_1\}.\]
	Clearly, $N$ is an ideal of $L_1 \times M$ and $N \subseteq \ker \tau_1 \phi_{L_1}$. Define $\gamma_1 : L_1 \to (L_1 \times M)/N$, $\gamma_1(l_1) = (l_1, {^M\{\mathcal{L}, G_1\}})N$ and $\gamma_2 : L_2 \to (L_1 \times M)/N$, $\gamma_2(l_2) = (l_1, (l_1, l_2) {^M\{\mathcal{L}, G_1\}})N$, where $\gamma \tau_1(l_1) = \tau_2(l_2)$. Clearly, $\ker \gamma_2 = 1$, let $l_1 \in \ker \gamma_1$, then $l_1 \in \ker \tau_1$ and $\beta(l_1) = 1$ this implies $l_1 = 1$, $\ker \gamma_1 = 1$. Thus $\gamma_1$ and $\gamma_2$ are monomorphisms. This conclude that the morphisms $(\gamma_1, I_{G_1}) : \tau_1 \to \overline{\tau_1 \phi_{L_1^\star}}$ and $(\gamma_2, \theta^{-1}) : \tau_2 \to \overline{\tau_1 \phi_{L_1}}$ are isoclinic monomorphisms.
\end{proof}
From the above propositions, we can easily obtained proof of the following theorem.
\begin{theorem}\label{main_thm}
	Let $\tau_1$ and $\tau_2$ be two relative Lie central extensions. Then the following are equivalent:
	\begin{enumerate}
		\item $\tau_1$ and $\tau_2$ are isoclinic.
		
		\item There exists a relative Lie central extension $\tau'$ together with isoclinic epimorphisms from $\tau'$ onto $\tau_1$ and $\tau_2$.
		
		\item There exists a relative Lie central extension $\tau''$ together with isoclinic monomorphisms from $\tau_1$ and $\tau_2$ into $\tau''$.
		
		\item There exists an abelian Lie algebra $M$, a relative Lie central extension $\tau'$, an isoclinic monomorphism from $\tau'$ into $\tau_1 \phi_{L_1} : L_1 \times M \to G_1$, and an isoclinic epimorphism from $\tau'$ onto $\tau_2$.
		
		\item There exists an abelian multiplicative Lie algebra with trivial bracket $M$, a relative Lie central extension $\tau''$, an isoclinic epimorphism from $\tau_1 \phi_{L_1} : L_1 \times M \to G_1$ onto $\tau''$, and an isoclinic monomorphism from $\tau_2$ into $\tau''$.
	\end{enumerate}
\end{theorem}

	\section{Schur multiplier of pair of multiplicative Lie algebras}
	
	In this section, we define the Schur multiplier of the pair of multiplicative Lie algebras and covering pair of the relative Lie central extension.
	

	Let $\xymatrix{1\ar[r] & R \ar[r]^{i} & F \ar[r]^{\beta} & G \ar[r] & 1}$ be the free presentation of the group $G$. Then the Schur multiplier of $G$, denoted by $M(G)$ is defined as
	\[M(G) \cong \dfrac{R \cap [F, F]}{[R, F]}.\]
	In \cite{Ellis_Pair}, Ellis defined the schur multiplier of the pair of groups $(G, N)$ as
	\[M(G, N) = \mathrm{Ker} (\mu : M(G) \to M(G/ N)),\]
	where $\mu$ is the natural epimorphism. Hence if $S$ is normal subgroup of $F$ such that $N \cong S/R$ then $M(G, N) \cong \dfrac{R \cap [S, F]}{[R, F]}.$
	
	In \cite{RLS}, Lal and Upadhyay, defined the Schur multiplier of multiplicative Lie algebras.
	\begin{definition}
		\cite{RLS} Let $E \equiv \xymatrix{1\ar[r] & R \ar[r]^{i} & F \ar[r]^{\beta} & G \ar[r] & 1}$ be free presentation of the multiplicative Lie algebra $G$. Then the Schur multiplier $\tilde{M}(G)$ of $G$ is defined as 
		\[\tilde{M}(G) \cong \dfrac{R \cap {^M[F, F]}}{^M[R, F]},\]
		where $^M[F, F] = (F \star F)[F, F]$ and $^M[R, F] = (R \star F)[R, F]$.
	\end{definition}
	
	\begin{theorem}\cite{RLS} \label{natural_exact_seq}
		To every Lie algebra extension $E \equiv \xymatrix{1\ar[r] & H \ar[r]^{i} & G \ar[r]^{\beta} & K \ar[r] & 1}$ there is an associated natural connecting homomorphism $\delta(E)$ from $\tilde{M}(K)$ to $\dfrac{H}{^M[H, G]}$
		and in turn the five term natural exact sequence
		$$ \tilde{M}(G) \xrightarrow{\makebox[1.2 cm]{$\tilde{M}(\beta)$}} \tilde{M}(K) \xrightarrow{\makebox[1cm]{$\delta(E)$}} \dfrac{H}{^M[H, G]} \xrightarrow{\makebox[0.8cm]{$\bar{i}$}} G_{ab} \xrightarrow{\makebox[0.8cm]{$\bar{\beta}$}} K_{ab} \parbox{0.8cm}{\rightarrowfill} 1.$$
	\end{theorem}

	The above theorem motivate us to define the Schur multiplier for the pair of multiplicative Lie algebras.
	\begin{definition}\label{schur_pair}
		Let $(H, G)$ be a pair of multiplicative Lie algebras. Then the Schur multiplier of $(H, G)$ is defined as 
		\[\tilde{M}(H, G) \cong \ker (\tilde{M}(\beta) : \tilde{M}(G) \to \tilde{M}(G/H)). \]
		If $S$ is an ideal of $F$ such that $H \cong S/R$, then 
		\[\tilde{M}(H, G)\cong \dfrac{R \cap {^M[S,F]}}{{^M[R,F]}}.\]
	\end{definition}
	
	
	\begin{proposition}\label{exact_sequence}
		Let $\xymatrix{1\ar[r] & R \ar[r] & F \ar[r] & G \ar[r] & 1}$ be a free presentation of the multiplicative Lie algebra $G$. If $S$ and $T$ are ideals of $F$ such that $S/R \cong H$ and $T/R \cong K$, where $T \subseteq S$, then the following sequence is exact
		$$1 \xrightarrow{ } \dfrac{R \cap {^M[T, F]}}{^M[R, F]} \xrightarrow{\makebox[0.4cm]{ }} \tilde{M}(H, G) \xrightarrow{\makebox[0.6cm]{$\beta$}} \tilde{M}\Big(\dfrac{H}{K}, \dfrac{G}{K}\Big) \xrightarrow{\makebox[0.6cm]{$\bar{i}$}} \dfrac{K \cap {^M[H, G]}}{^M[K, G]} \xrightarrow{ }  1.$$
	\end{proposition}
	\begin{proof}
		Definition \ref{schur_pair}, gives
		$$\tilde{M}(H, G) \cong \dfrac{R \cap {^M[S, F]}}{^M[R, F]},~~ \tilde{M}\Big(\dfrac{H}{K}, \dfrac{G}{K}\Big) \cong \dfrac{T \cap {^M[S, F]}}{^M[S, F]}$$
		and 
		$$ \dfrac{K \cap {^M[H, G]}}{^M[K, G]} \cong \dfrac{T \cap {^M[S, F]} R}{^M[T, F] R}.$$
		Thus, the exactness of the sequence with natural homomorphisms can be easily verified
		\begin{multline*}
			1 \xrightarrow{\makebox[0.6cm]{ }} \dfrac{R \cap {^M[T, F]}}{^M[R, F]} \xrightarrow{\makebox[0.4cm]{ }} \dfrac{R \cap {^M[S, F]}}{^M[R, F]} \xrightarrow{\makebox[0.6cm]{$\beta$}} \dfrac{T \cap {^M[S, F]}}{^M[S, F]}
			\xrightarrow{\makebox[0.6cm]{ }} \dfrac{T \cap {^M[S, F]} R}{^M[T, F]R} \xrightarrow{ \makebox[0.6cm]{ }}  1.	\end{multline*}
	\end{proof}
	
	\begin{corollary}
		Let $H$ be an ideal of $G$ such that $\tilde{M}(H, G)= 1$. Then $\tilde{M}\Big(\dfrac{H}{K}, \dfrac{G}{K}\Big) \cong \dfrac{K \cap {^M[H, G]}}{^M[K, G]}.$ 
	\end{corollary}
	
	\begin{corollary}
		Let  $ \xymatrix{1\ar[r] & I \ar[r] & G \ar[r] & K \ar[r] & 1}$ be central extension of $K$ such that $I \subseteq {^M[H, G]}$ and $\tilde{M}\Big(\dfrac{H}{K}, \dfrac{G}{K}\Big) = 1$. Then $G$ and $K$ are isomorphic.
	\end{corollary}
	
	\begin{definition}\label{covering_pair}
		Let $\tau : L \to G$ be a relative Lie central extension of the pair $(H, G)$. Then $\tau : L  \to G$ is called a multiplicative covering pair of $(H, G)$, if there exists a ideal $I$ of $L$ such that 
		\begin{enumerate}
			\item $I \subseteq \bar{\mathcal{Z}}(G, L) \cap {^M\{L, G\}}$
			\item $I \cong \tilde{M}(H, G)$ and $H \cong L /I$.
		\end{enumerate}
		
	\end{definition}


	Let $(H, G)$ be a pair of multiplicative Lie algebras. Let $$E \equiv \xymatrix{1\ar[r] & R \ar[r]^{i} & F \ar[r]^{\beta} & G \ar[r] & 1}$$ be a free presentation of $G$ and $S$ be an ideal of $F$ such that $H \cong S/R$. Then the mapping $\tilde{\tau} : S/{^M[R,F]} \to F/R$ given by $\tau (s^M[R,F]) = sR$ is a relative Lie central extension of the pair $(H, G)$ with the given actions
	\[ \dfrac{F}{R} \times \dfrac{S}{^M[R, F]} \to \dfrac{S}{^M[R, F]}  \] 
	\[(fR, s(^M[R,F])) \mapsto s(^M[R,F])^{fR} =  s^f ({^M[R,F]}) \]
	\[\langle fR, s(^M[R,F]) \rangle \mapsto fR \star s(^M[R,F])  =  (f \star s)^M[R,F]\]
	and 
	\[ \dfrac{S}{^M[R,F]} \times \dfrac{F}{R} \to \dfrac{F}{R}  \] 
	\[(s({^M[R,F]}), fR) \mapsto fR^{s(^M[R,F])} =  f^s R\]
	\[\langle s(^M[R,F], fR) \rangle \mapsto s(^M[R,F]) \star {fR} =  (s\star f)R.\]

	\begin{theorem} \label{commutative_diagram}
		Let	$E \equiv \xymatrix{1\ar[r] & R \ar[r]^{i} & F \ar[r]^{\beta} & G \ar[r] & 1}$ be a free presentation of $G$ and $H$ be an ideal of $G$ such that $H \cong S/R$, where $S$ is an ideal of $F$. If $\delta : L \to G$ is a relative Lie central extension of the pair $(H, G)$, then there exist a homomorphism $\phi$ from $S/[R,F](S \star F)$ to $L$ such that the following diagram is commutative:
		\[
		\xymatrix{
			1 \ar[r] & \dfrac{R}{^M[R, F]} \ar[r]^{ } \ar[d]^{\phi|} & \dfrac{S}{^M[R, F]} \ar[r]^{\hspace{1 cm}\tilde{\tau}} \ar[d]^{\phi} & H \ar[r] \ar[d]^{ } & 1 \\
			1 \ar[r] & \mathrm{Ker}(\delta) \ar[r]^{ } & L \ar[r]^{\delta} & H \ar[r] & 1,
		}
		\]
		where $\tilde{\tau}$ is the relative Lie central extension, described as above. Moreover, if $L$ is a perfect multiplicative Lie algebra and $\tilde{{\Phi}}(L) \neq L$, then $\phi$ is a surjective homomorphism.
	\end{theorem}
	\begin{proof}
		Define $p : G \to H$ by $p(g) = g$ if $g \in H$, otherwise $p(g) = 1$. Then $p$ induce a homomorphism $\chi : F \to H$. Since $\delta (L) = H$, by lifting property of multiplicative Lie algebras, there exist a homomorphism $\psi : F \to L$ such that $\delta \psi = \chi$. Now, restriction of $\psi$ and $\chi$ on S, gives the following commutative diagram
		\[
		\xymatrix{
			S \ar[rd]^{\hspace{.1 cm}\chi|} \ar[d]_{\psi|}  \\
			L \ar[r]^{ \delta} & L \ar[r] & 1
		}
		\]
		Let $r \in R$, then $\psi(r) \in \mathrm{Ker}(\delta) \leqslant \bar{\mathcal{Z}}(L, G) \leqslant \mathcal{Z}(L)$. For $x \in F$, we have $\psi ([x, r]) = [\psi(x), \psi(r)] = 1$ and $\psi(x \star r) = \psi(x) \star \psi(r) = 1$, $^M[R, F] \leqslant \mathrm{Ker}(\psi|)$, this induces a map $\phi : \dfrac{S}{^M[R, F]} \to L$. It can be easily seen that $\phi\Big(\dfrac{R}{^M[R, F]}\Big) \leqslant \mathrm{Ker}(\delta)$. Therefore $\phi$ can be restrict to $\dfrac{R}{^M[R,F]}$. 
		
		Let $l \in L$, then $\delta(l) = \tilde{\tau}(\bar{s}) = \delta \phi(\bar{s})$ for some $\bar{s} \in S/^M[R, F]$, this implies $l = \phi(\bar{s}) t$ for some $t \in \mathrm{Ker}(\delta)$. Thus
		\[ L = \big\langle \phi(\bar{s}), t : \bar{s} \in \dfrac{S}{^M[R, F]}, t \in \mathrm{Ker}(\delta) \big\rangle\] 
		Since $L = {^M[ L, L]} $, gives 
		\[ \mathrm{Ker}(\delta) \leqslant \bar{\mathcal{Z}}(L, G) \cap {^M[L,L]} \leqslant \mathcal{Z}(L) \cap {^M[L, L]} \leqslant \tilde{{\Phi}}(L).\]
		Thus, by Theorem \ref{Frattini_Non}
		\[ L = \big\langle \phi(\bar{s}): \bar{s} \in \dfrac{S}{^M[R, F]} \big\rangle,\] 
		i.e. $L$ is generated by $\im \phi $. Therefore, $\phi$ is a surjective homomorphism.
	\end{proof}
	
	
	\begin{corollary}
		Let $\tau : L \to G$ be a multiplicative covering pair of $(H, G)$ such that $L$ is perfect and $\tilde{{\Phi}}(L) \neq L$. Then $L$ is a homomorphic image of $\dfrac{S}{^M[R, F]}$.
	\end{corollary}
	

	\begin{proposition}
		Let $(H, G)$ be a pair of multiplicative Lie algebras and \\ $E \equiv \xymatrix{1\ar[r] & R \ar[r]^{i} & F \ar[r]^{\beta} & G \ar[r] & 1}$ be a free presentation of $G$ and $S$ is an ideal of $F$ such that $H \cong S/R$. If $T$ is a normal subgroup of $F$ such that
		\[\dfrac{R}{^M[R, F]} = \dfrac{R \cap {^M[S,F]}}{^M[R, F]} \times \dfrac{T}{^M[R,F]}\] 
		then the map $\tau : S/T \to F/R$ given by $\tau (sT) = sR$ is a covering pair of $(H,G)$.
	\end{proposition}
	\begin{proof}
		Firstly, we will show that $F/R$ is acting on $S/T$. Let group actions of $F/R$ on $S/T$ and $S/T$ on $F/R$ be given by $(fR, sT) \mapsto {^{fR}}sT = {^fs}T$ and $(s'T, f'R) \mapsto {^{s'T}}f'R = {^{s'}f'}R$, resepectively. Where as Lie bracket action is defined as 
		\[ \langle ~, ~ \rangle : \dfrac{F}{R} \times \dfrac{S}{T} \to \dfrac{S}{T}; ~ \langle fR, sT \rangle \mapsto (f \star s) T.\]
		Then it is easy to check that with the above action, $F/R$ is acting on $S/T$. Now we will show that $\tau$ is a relative Lie central extension of the pair $(H,G)$.
		\[\tau(^{fR}sT) = \tau (^fsT) = {^fs}R = fR sR f^{-1}R = fR \tau(sT) f^{-1}R\]
		\[\tau (\langle fR, sT \rangle) = \tau((f\star s)T) = (f\star s)R = fR \star sR = fR \star \tau(sT)\]
		\[^{\tau(sT)}s'T= {^{sR}}s'T= {^ss'} T = sT s'T s^{-1}T\]
		\[\langle \tau(sT), s'T \rangle = \langle sR, s'T \rangle = (s\star s')T = sT \star s'T\]
		\[^{\tau(sT)}fR = {^{sR}}fR = {^sf}R = {^{sT}}fR.\]
		Clearly, $\tau(S/T)= S/R$ and $F/R$ acts trivially on $\ker \tau$.
		
		\noindent Let $L^\star = S/T$ and $I = R/T$, then $H \cong S/R$. Then 
		\[I = \dfrac{R}{T} \cong \dfrac{R \cap {^M[S, F]}}{^M[R, F]} \cong \mathcal{Z}(G,H).\]
		Clearly, $I \subseteq \mathcal{Z}(G, H)$ for $F/R$ acts trivially on $\ker (\tau)$.
		\begin{align*}
		I = \dfrac{R}{T} \subseteq \dfrac{^M[S,F]T}{T} &  
			 = \big\langle {^fs}T (sT)^{-1}(f'T\star s'T) : f,f' \in F~ s,s' \in S \big\rangle \\
			& = \big\langle {^{fR}}sT (sT)^{-1}  (f' \star s')T  : f,f' \in F~ s,s' \in S \big\rangle \\
			& = \big\langle {^{fR}}sT (sT)^{-1} \langle f'R, s'T \rangle : f,f' \in F~ s,s' \in S \big\rangle.\\
			& = {^M\{G, L^\star\}}.
		\end{align*}
		Thus $I \subseteq \mathcal{Z}(G, H) \cap {^M\{G, L^\star\}}$. Therefore, $\tau : S/T \to F/R$ is the multiplicative covering pair of $(H,G)$.
	\end{proof}
	
	
	\begin{proposition} \label{perfect_L}
		Let $(H, G)$ be a perfect pair of multiplicative Lie algebras and $\tau : L \to G$ be multiplicative covering pair. Then $^M\{L, G\} = L$.
	\end{proposition}
	\begin{proof}
		Given $\tau : L \to G$ is a multiplicative covering pair, this implies $L/ \mathrm{Ker}(\tau) \cong H$ and $\mathrm{Ker}(\tau) \subseteq {^M\{L, G\}}$. Thus, isomorphism theorem gives
		\[\dfrac{L/\mathrm{Ker}(\tau)}{^M\{L, G\}/\mathrm{Ker}(\tau)} \cong \dfrac{H}{^M\{H, G\}}.\]
		Since, $H = {^M\{H, G\}}$, we have $L = {^M\{L, G\}}$.
	\end{proof}
	

\end{document}